\documentclass[11pt,a4paper]{amsart}

\usepackage[utf8x]{inputenc} 
\usepackage{enumerate}  
\usepackage{fancyvrb}    
\usepackage{color,colortbl}
\usepackage[dvipsnames]{xcolor}
\usepackage{hyperref}  
\usepackage{verbatim}
                 \hypersetup{ pdfborder={0 0 0}, 
                              colorlinks=true, 
                              linkcolor=red, 
                              urlcolor=red,
                              citecolor=red,
                              linktoc=page,
                              pdfauthor={Giovanni Staglian{\`o}}, 
                              pdftitle={A package for computations with classical resultants}
                            }                         
\usepackage{mathptmx}                      
   
\theoremstyle{plain} 
\newtheorem{proposition}{Proposition}[section] 
\newtheorem{theorem}[proposition]{Theorem}
\theoremstyle{definition} 
\newtheorem{definition}[proposition]{Definition} 
\newtheorem{example}[proposition]{Example} 

\theoremstyle{remark} 
\newtheorem{remark}[proposition]{Remark}

\newcommand{\GG}{{\mathbb{G}}} 
\newcommand{\ZZ}{{\mathbb{Z}}} 
\newcommand{\CC}{{\mathbb{C}}}  
\newcommand{\Aa}{{\mathbb{A}}} 
         
\newcommand{\PP}{{\mathbb{P}}}          

\providecommand{\Res}{\mathop{\rm Res}\nolimits} 
\providecommand{\Dis}{\mathop{\rm Disc}\nolimits} 
\providecommand{\det}{\mathop{\rm det}\nolimits} 

\numberwithin{equation}{section}

\title{A package for computations with classical resultants} 
\author[G. Staglian\`o]{Giovanni Staglian\`o} 
\address{Dipartimento di Ingegneria Industriale e Scienze Matematiche, Universit\`a Politecnica delle Marche, via Brecce Bianche, I-60131 Ancona, Italy}
\date{\today} 
\email{\href{mailto:giovannistagliano@gmail.com}{giovannistagliano@gmail.com}} 

\begin{document}

\begin{abstract} 
We present the \emph{Macaulay2} package 
\href{https://github.com/Macaulay2/M2/blob/master/M2/Macaulay2/packages/Resultants.m2}{\emph{Resultants}},
which provides 
 commands
for the effective computation of multivariate resultants, discriminants, and 
Chow forms. 
We provide some background for the algorithms implemented 
and show, with a few examples, 
 how the package works.
\end{abstract}

\maketitle

\section*{Introduction}
The \emph{resultant} characterizes the existence of nontrivial solutions 
for a square system of homogeneous polynomial equations
as a condition on the coefficients.  
One of its important features is 
 that it
can be used to compute elimination ideals and to solve polynomial equations. Indeed,
it provides one of the two main 
tools in elimination theory,
along with Gr\"obner bases. 
The resultant of the system of equations given 
by the partial derivatives of a complex homogeneous polynomial $F$ 
is  called (up to a constant factor)
the \emph{discriminant} of $F$. 
It 
characterizes the existence of singular points 
in the projective hypersurface $V(F)$
as a condition on the coefficients of $F$.
In this special case,
all polynomial equations have the same total degree.
Every time 
the system of equations consists of
  $n+1$ polynomial equations 
of the same total degree $d$,
the resultant  has a further interesting 
                       property: 
it can be expressed as a polynomial of degree $d^n$ 
in the $(n+1)\times(n+1)$ minors 
of an $(n+1)\times\binom{n+d}{n}$ matrix, the coefficient matrix of the 
system of equations. 
This allows us to write down a \emph{generic} resultant in a more compact form. 
The polynomial of degree $d^n$ so obtained 
is geometrically interpreted as 
the \emph{Chow form} of the $d$-th Veronese embedding of $\PP^n$.

The package \href{https://github.com/Macaulay2/M2/blob/master/M2/Macaulay2/packages/Resultants.m2}{\emph{Resultants}}, 
included with \emph{Macaulay2}  \cite{macaulay2},
provides commands for the explicit computation of resultants and discriminants.
The main algorithm used is based on the so-called \emph{Poisson formula}, 
which reduces the computation of the resultant of $n+1$  equations to the 
product of the resultant of $n$ equations with the determinant 
of an appropriate matrix.
This algorithm requires a certain genericity condition on the input polynomials,
achievable with a generic change of coordinates.
The package  also 
includes  tools for working with
Chow forms and more generally with \emph{tangential Chow forms}.

In Section~\ref{section1}, from a more computational point of view,
we give some background information on the general theory 
of 
resultants, discriminants, and Chow forms.
In Section~\ref{section2}, 
we briefly illustrate how to use the package 
with the help of some examples;
    more detailed information and examples 
can be found in its documentation.

\section{Overview of classical resultants}\label{section1}
In this section, we present an overview of some classically well known facts on the theory 
of resultants for forms in several variables.
For details and proofs,
we refer mainly to 
\cite{gelfandbook} and \cite{coxbook}; other references are
\cite{jouanolou1,jouanolou2,Waerden1950,demazure1,EMIRIS19993,bgw-amer-88,Bus2014}, and \cite{coxbook2} 
for the 
case of two bivariate polynomials.
\subsection{Resultants}
Suppose we are given $n+1$ homogeneous polynomials $F_0,\ldots,F_n$ 
 in $n+1$ variables $x_0,\ldots,x_n$ over the complex field $\CC$. 
 For $i=0,\ldots,n$, 
 let $d_i$ denote the total degree of $F_i$ so that we can write 
 $F_i = \sum_{|\alpha|=d_i}c_{i,\alpha} x^{\alpha}$, where  
 $x^{\alpha}$ denotes $x_0^{\alpha_0}\cdots x_n^{\alpha_n}$.
 For each pair of indices $i,\alpha$, we introduce a variable $u_{i,\alpha}$ and form 
 the \emph{universal ring of coefficients} $\mathbb{U}_{d_0,\ldots,d_n}:=\ZZ[u_{i,\alpha}:i=0,\ldots,n,|\alpha|=d_i]$.
 If $P\in\mathbb{U}_{d_0,\ldots,d_n}$,
 we denote by $P(F_0,\ldots,F_n)$ 
 the element in $\CC$ obtained by replacing each variable $u_{i,\alpha}$ 
 with the corresponding coefficient $c_{i,\alpha}$.
\begin{theorem}[\cite{gelfandbook,coxbook}]
 If we fix positive degrees $d_0,\ldots,d_n$, then there is a unique 
 polynomial $\Res=\Res_{d_0,\ldots,d_n}\in\mathbb{U}_{d_0,\ldots,d_n}$ which has the following properties:
 \begin{enumerate}
  \item if $F_0,\ldots,F_n\in\CC[x_0,\ldots,x_n]$ are homogeneous of degrees $d_0,\ldots,d_n$,
  then the equations 
  \begin{equation*}
   F_0 = 0,\ldots,F_n = 0
  \end{equation*}
 have a nontrivial solution over $\CC$ (i.e., $\emptyset\neq V(F_0,\ldots,F_n)\subset\PP^n_{\CC}$) if and only if $\Res(F_0,\ldots,F_n)=0$;
  \item $\Res$ is irreducible, even when regarded as a polynomial over $\CC$;
  \item $\Res(x_0^{d_0},\ldots,x_n^{d_n})=1$.
 \end{enumerate}
\end{theorem}
\begin{definition} We call $\Res(F_0,\ldots,F_n)$ the \emph{resultant} of $F_0,\ldots,F_n$. \end{definition}
\begin{remark}\label{remark-specialization}
 If $A$ is any commutative ring, we define the resultant of $n+1$ 
 homogeneous polynomials $F_0,\ldots,F_n\in A[x_0,\ldots,x_n]$ again as 
 $\Res(F_0,\ldots,F_n)\in A$, i.e., by specializing the coefficients of the integer polynomial $\Res$.
 Thus, 
 the formation of resultants
 commutes with specialization.
\end{remark}
\begin{example}
The resultant is a direct generalization of the determinant. Indeed,
 if $d_0=\cdots=d_n=1$, then $\Res(F_0,\ldots,F_n)$ equals the determinant of 
 the $(n+1)\times(n+1)$ coefficient matrix.
\end{example}
\begin{proposition}[\cite{jouanolou1,jouanolou2}]\label{propertiesRes} The following hold:
 \begin{enumerate}
  \item\emph{(Homogeneity)} For a fixed $j$ between $0$ and $n$, $\Res$ is homogeneous 
  in the variables $u_{j,\alpha}$, $|\alpha|=d_j$, of degree $d_0\cdots d_{j-1}d_{j+1}\cdots d_n$;
  hence its total degree  is $\sum_{j=0}^n d_0\cdots d_{j-1}d_{j+1}\cdots d_n$.
  \item\emph{(Simmetry)} If $\sigma$ is a permutation of $\{0,\ldots,n\}$, then 
         \[\Res(F_{\sigma(0)},\ldots,F_{\sigma(n)}) 
         = \mathrm{sign}(\sigma)^{d_0\cdots d_n}\,\Res(F_0,\ldots,F_n) . \]
  \item\emph{(Multiplicativity)} If $F_j=F_j'F_j''$, then 
  we have 
  \[\Res(F_0,\ldots,F_j,\ldots,F_n) = 
  \Res(F_0,\ldots,F_j',\ldots,F_n) \, 
  \Res(F_0,\ldots,F_j'',\ldots,F_n) . \]
  \item\emph{($\mathrm{SL}(n+1)$-invariance)}\label{relativeInv}  For each $(n+1)\times(n+1)$ matrix $A$ over $\CC$, 
  we have \[\Res(F_0(A\,x),\ldots,F_n(A\,x)) 
           = {\det(A)}^{d_0\cdots d_n}\Res(F_0(x),\ldots,F_n(x)) , \]
           where $A\,x$ denotes the product of $A$ with the column vector $(x_0,\ldots,x_n)^t$.
                    \item \emph{(Elementary transf.)} If $H_i$ is homogeneous  
  of degree $d_j - d_i$, then
                    \[ \Res(F_0,\ldots,F_j+\sum_{i\neq j} H_i F_i,\ldots,F_n)
  =
  \Res(F_0,\ldots,F_j,\ldots,F_n) . \] 
  \end{enumerate}
\end{proposition}
\begin{remark}
On the product $\Aa^M\times \PP^n=\mathrm{Spec}(\CC[u_{i,\alpha}])\times\mathrm{Proj}(\CC[x_0,\ldots,x_n])$, 
      where $M=\sum_{i=0}^n\binom{n+d_i}{n}$, we have an  \emph{incidence variety}
\begin{equation*}
  W:=\{((c_{i,\alpha}),p)\in\Aa^M\times \PP^n : p\in V(\sum_{|\alpha|=d_0}c_{0,\alpha}x^{\alpha},\ldots,\sum_{|\alpha|=d_n}c_{n,\alpha}x^{\alpha})\}. 
\end{equation*}
 The first projection  $\pi_1:W\rightarrow\Aa^M$ is 
 birational onto its image,  whereas 
 all the fibers of the second projection $\pi_2:W\rightarrow\PP^n$ 
 are linear subspaces of dimension $M-n-1$. It follows that 
 $W$ is a smooth irreducible variety 
 which is  birational to
 $\overline{\pi_1(W)}=\pi_1(W)=V(\Res_{d_0,\ldots,d_n})\subset\Aa^M$.
\end{remark}
The following result is called \emph{Poisson formula} and allows one to compute resultants 
inductively.
\begin{theorem}[\cite{jouanolou1}; see also \cite{coxbook}]\label{PoissonFormula}
Let $f_i(x_0,\ldots,x_{n-1}):=F_i(x_0,\ldots,x_{n-1},1)$ and 
 $\overline{F_i}(x_0,\ldots,x_{n-1}):=F_i(x_0,\ldots,x_{n-1},0)$.
 If $\Res(\overline{F_0},\ldots,\overline{F_{n-1}})\neq 0$, 
 then the quotient ring $A=\CC[x_0,\ldots,x_{n-1}]/(f_0,\ldots,f_{n-1})$ 
 has dimension $d_0\cdots d_{n-1}$ as a vector space over $\CC$, and 
 \begin{equation}\label{PoissonForEq}
  \Res(F_0,\ldots,F_n)
   = \Res(\overline{F_0},\ldots,\overline{F_{n-1}})^{d_n}\,\det(m_{f_n}:A\rightarrow A) , 
   \end{equation}
where $m_{f_n}:A\rightarrow A$ is the linear map given by multiplication by $f_n$.
\end{theorem}
With the same hypotheses as Theorem \ref{PoissonFormula}, 
a monomial basis for $A$ over $\CC$ (useful in the implementation)
can be 
constructed as explained
 in \cite[Chapter~2, \S~2]{coxbook}.
Note also that 
we have
\begin{equation}
 \det(m_{f_n}:A\rightarrow A)=\prod_{p\in V} f_n(p)^{\mathrm{mult}_p(V)} 
,\end{equation} where $V=V(f_0,\ldots,f_{n-1})$. 
 
We now describe 
the most popular way to compute resultants, which is
due to Macaulay \cite{macaulaypaper1902}. 
Let $\delta=\sum_{i=0}^n d_i - n$ and $N=\binom{n+\delta}{n}$.  
We can divide the monomials $x^{\alpha}$ of total degree $\delta$  into the 
$n+1$ mutually disjoint sets 
\[ S_i:=\{x^{\alpha} : |\alpha|=\delta,\ \mathrm{min}\{j:x_j^{d_j}|x^{\alpha}\}=i\},  \mbox{ for }i=0,\ldots n.\] A monomial 
$x^{\alpha}$ of total degree $\delta$ is called \emph{reduced} if $x_i^{d_i}$ divides $x^{\alpha}$ 
for exactly one~$i$. Consider the following $N$ homogeneous polynomials of degree $\delta$:
\begin{eqnarray}\label{eqs4.1}
 x^{\alpha}/x_i^{d_i}\,F_i, \mbox{ for }i=0,\ldots n\mbox{  and  }x^{\alpha}\in S_i.
\end{eqnarray}
By regarding the monomials of total degree $\delta$ as unknowns, the polynomials 
in \eqref{eqs4.1} form 
a system of $N$ linear equations in $N$ unknowns. 
Let $\mathbb{D}=\mathbb{D}(F_0,\ldots,F_n)$ denote
the coefficient matrix of this linear system, and
let $\mathbb{D}'(F_0,\ldots,F_n)$ denote the submatrix of $\mathbb{D}$ 
obtained by deleting all rows and columns corresponding to reduced monomials.
The following result is called the \emph{Macaulay formula} 
and allows one to compute the resultant as a quotient of two determinants.
\begin{theorem}[\cite{macaulaypaper1902}; see also \cite{jouanolou2,coxbook}]\label{MacaulayFormula}
The following formula holds: 
 \begin{equation}\label{macaulayEqFor}
\det(\mathbb{D}(F_0,\ldots,F_n)) = \Res(F_0,\ldots,F_n)\,\det(\mathbb{D}'(F_0,\ldots,F_n)) .
   \end{equation}
\end{theorem}
In several special cases,
the resultant can be expressed as a single determinant (see 
\cite[Chapter~13, Proposition~1.6]{gelfandbook}). 
We also mention that besides \eqref{macaulayEqFor},
there are other ways to represent resultants as quotients:
these include \emph{Bezoutians} \cite{elkadi:inria-00073109}
and \emph{Dixon matrices} \cite{KaSaYa94}; see also \cite{EMIRIS19993} and \cite[p.~110]{coxbook}.
However, 
all these matrices are usually of much larger size than those involved by the 
Poisson formula \eqref{PoissonForEq}, as shown in the following simple example (see 
 \cite{EMIRIS19993},
for a comparison between Macaulay and other resultant matrices).
\begin{example}
Let 
$F_0=x^{3}+y^{2} z$, $F_1=x y+y^{2}+x z+y z$, $F_2=y^{4}+z^{4}\in\CC[x,y,z]$.
The Poisson formula expresses $\Res(F_0,F_1,F_2)$ as the following product of determinants:
\[
\Res(F_0,F_1,F_2) = 
\left( 
           1^2\,  \det  \bgroup\begin{pmatrix}1&
     0&
     0\\
     1&
     1&
     0\\
     0&
     1&
     1\\
     \end{pmatrix}\egroup
 \right)^{4}\,
 \det \bgroup\begin{pmatrix}1&
     0&
     0&
     0&
     0&
     0\\
     0&
     2&
     0&
     0&
     0&
     0\\
     0&
     0&
     2&
     0&
     0&
     0\\
     0&
     1&
     0&
     1&
     0&
     0\\
     1&
     0&
     0&
     0&
     2&
     0\\
     0&
     {-1}&
     0&
     1&
     0&
     2\\
     \end{pmatrix}\egroup = 16 . 
 \]
 The Macaulay formula yields the same result as a quotient $\det(\mathbb{D})/\det(\mathbb{D}')$,
 where 
 $\mathbb{D}$ and $\mathbb{D}'$ are square matrices of size $36\times 36$ and $10\times 10$,
 respectively.
\end{example}
\subsection{Discriminants}
Let $F=\sum_{|\alpha|=d}c_{\alpha} x^{\alpha}\in\CC[x_0,\ldots,x_n]$ 
be a homogeneous polynomial of a certain
degree $d$.
As above, for each index $\alpha$ we introduce a variable $u_{\alpha}$
and form the universal ring of coefficients 
 $\mathbb{U}_d:=\CC[u_{\alpha}:|\alpha|=d]$. Then 
 one can show that,
 up to sign,
 there is a unique polynomial
 $\Dis=\Dis_{d}\in\mathbb{U}_d$ 
 which has the following properties:
 \begin{enumerate}
  \item if $F\in\CC[x_0,\ldots,x_n]$ is homogeneous of degrees $d$,
  then the equations 
  \[
   \partial F/\partial x_0 = 0,\ldots,\partial F/\partial x_n = 0
  \]
 have a nontrivial solution over $\CC$ 
 (i.e., the hypersurface defined by $F$ is singular) if and only if $\Dis(F)=0$;
  \item $\Dis$ is irreducible, even when regarded as a polynomial over $\CC$.
 \end{enumerate}
\begin{proposition}[\cite{gelfandbook}]\label{discriProp} 
Up to sign, we have the formula 
\begin{equation}\label{discrForm}
 \Dis(F) 
 = c_{d,n}\, \Res(\frac{\partial F}{\partial x_0},\ldots,\frac{\partial F}{\partial x_n}),
 \mbox{ where }
  c_{d,n} = d^{\frac{(-1)^{n+1} - (d-1)^{n+1}}{d}} .
\end{equation}
\end{proposition} 
\begin{definition} We call the polynomial defined by \eqref{discrForm} the \emph{discriminant} of $F$. \end{definition}
\begin{proposition}[\cite{gelfandbook}] 
            The following hold:
 \begin{enumerate}
  \item The polynomial  $\Dis$ is homogeneous  
  of degree $(n+1)(d-1)^n$.
  \item For each   $(n+1)\times(n+1)$ matrix $A$ over $\CC$, 
  we have \[\Dis(F(A\,x)) = {\det(A)}^{d(d-1)^n}\Dis(F(x)), \]
           where $A\,x$ denotes the product of $A$ with the column vector $(x_0,\ldots,x_n)^t$.
 \end{enumerate}
\end{proposition}
 Geometrically, we have the following interpretation.
\begin{proposition}[\cite{gelfandbook}]
 The discriminant hypersurface  $V(\Dis_d)$ 
 in the space of forms of degree $d$ on $\PP^n$
  coincides with 
 the dual variety of the $d$-th Veronese embedding of $\PP^n$. 
 \end{proposition}
 
\subsection{Chow forms}
Let $X\subset \PP^n$ be an irreducible subvariety of dimension $k$ and degree $d$. Consider 
the subvariety ${Z}(X)$ 
in the Grassmannian $\GG(n-k-1,\PP^n)$ 
of all $(n-k-1)$-dimensional projective subspaces of $\PP^n$ that intersect $X$.
It turns out that ${Z}(X)$ is an irreducible hypersurface of degree $d$;
 thus ${Z}(X)$ is defined by the vanishing 
of some element $R_X$, unique up to a constant factor, in the homogeneous component 
of degree $d$ of the coordinate ring of the Grassmannian $\GG(n-k-1,\PP^n)$ in the Pl\"{u}cker embedding. 
This element is called the \emph{Chow form} of  $X$. It is notable that $X$ can 
be recovered from its Chow form. See \cite[Chapter~3, \S~2]{gelfandbook} for details.

Consider the product $\PP^k\times X$ as a subvariety of $\PP^{(k+1)(n+1)-1}$
via the Segre embedding. Identify $\PP^{(k+1)(n+1)-1}$ with the 
projectivization 
$\PP(\mathrm{Mat}(k+1,n+1))$ 
 of the space of $(k+1)\times(n+1)$ matrices 
  and consider the natural projection 
  $\rho:\PP(\mathrm{Mat}(k+1,n+1))\dashrightarrow \GG(k,n)\simeq \GG(n-k-1,n)$.
  The following result is called the \emph{Cayley trick}. 
  \begin{theorem}[\cite{gelfandbook}; see also \cite{Weyman1994}]\label{CayleyTrickTeo}
   The dual variety of $\PP^k\times X$ coincides  with 
   the closure $\overline{\rho^{-1}(Z(X))}$, where $Z(X)\subset \GG(n-k-1,n)$ 
   is the hypersurface defined by the Chow form of $X$.
  \end{theorem}
The defining polynomial of the hypersurface 
$\overline{\rho^{-1}(Z(X))}\subset\PP(\mathrm{Mat}(k+1,n+1))$ 
is called \emph{$X$-resultant};
it provides another way of writing the Chow form of $X$.

Now, let $F_0,\ldots,F_n$ be $n+1$ generic homogeneous polynomials 
on $\PP^n$ of the same degree $d>0$,
and let $\mathbb{M}=\mathbb{M}(F_0,\ldots,F_n)$ be 
the $(n+1)\times N$ matrix  
of the coefficients 
of these polynomials,  $N=\binom{n+d}{n}$.
We consider the projection 
$\rho_{n,d}:\PP(\mathrm{Mat}(n+1,N))\dashrightarrow\GG(n,N-1)\simeq \GG(N-n-2,N-1)$ 
defined by the maximal minors of $\mathbb{M}$.
\begin{proposition}[\cite{gelfandbook,coxbook}]\label{XresVeronese}
 The hypersurface of degree $(n+1)d^n$ in $\PP(\mathrm{Mat}(n+1,N))$ 
 defined by the resultant $\Res(F_0,\ldots,F_n)$ coincides 
 with the closure $\overline{\rho_{n,d}^{-1}(V(R_{n,d}))}$, where $R_{n,d}$ denotes 
 the Chow form of the $d$-th Veronese embedding of $\PP^n$.
 In particular, $\Res(F_0,\ldots,F_n)$ is a polynomial 
 in the maximal minors of $\mathbb{M}$.
\end{proposition}

\section{Implementation}\label{section2}
In this section, we illustrate briefly some of the methods available 
in the package 
\href{https://github.com/Macaulay2/M2/blob/master/M2/Macaulay2/packages/Resultants.m2}{\emph{Resultants}}, 
included with \emph{Macaulay2}  \cite{macaulay2}. 
We refer to the package documentation (which can be viewed with \texttt{viewHelp Resultants}) 
for more 
details and examples. 

One of the main methods is 
 \texttt{resultant}, which accepts as input a list of $n+1$ homogeneous polynomials 
 in $n+1$ variables with coefficients in some commutative ring $A$ 
 and returns an element of $A$, the resultant of the polynomials.
 There are no limitations on the ring $A$  because of Remark~\ref{remark-specialization}.
 The algorithms implemented are 
the Poisson formula (Theorem \ref{PoissonFormula}) and 
the Macaulay formula (Theorem \ref{MacaulayFormula}). 
The former 
 is used by default since  it is typically 
                                         faster, 
  while for the latter one has to set the Algorithm option:
  \texttt{resultant(...,Algorithm=>"Macaulay")}. 
The method can also be configured to involve 
interpolation of multivariate polynomials (see \cite{CannyManocha}),
i.e., it can reconstruct the polynomial resultant from its values at a sufficiently 
large number of points, which in turn 
are evaluated using the same formulas.
The main derived method  is 
\texttt{discriminant}, which applies the formula \eqref{discrForm} to compute discriminants of homogeneous polynomials. 
\begin{example}
In the following code, we take two forms $F,G$ of degree $6$ on $\PP^3$. 
We first verify that $\Dis(F)=0$ and  $\Dis(G)\neq 0$ and then we 
compute 
the intersection 
of the pencil generated by $F$ and $G$ 
with the discriminant hypersurface in the space of forms of degree $6$ on $\PP^3$,
which is a hypersurface of degree $500$ in $\PP^{83}$. 
(The algorithm behind these calculations
is the Poisson formula; this is one of the cases where 
the Macaulay formula is much slower).
{\footnotesize 
\begin{Verbatim}[commandchars=&\{\}]
Macaulay2, version 1.10
with packages: ConwayPolynomials, Elimination, IntegralClosure, InverseSystems,
               LLLBases, PrimaryDecomposition, ReesAlgebra, TangentCone
i1 : loadPackage "&colore{blue}{Resultants}";
i2 : ZZ[w,x,y,z]; (F,G) = (w^6+x^6+y^6+w*x*y^4,w^6+x^6+y^6+z^6)
       6    6        4    6   6    6    6    6
o3 = (w  + x  + w*x*y  + y , w  + x  + y  + z )
o3 : Sequence
i4 : time &colore{blue}{discriminant} F
     &colore{Sepia}{-- used 0.0179806 seconds}
o4 = 0
i5 : time &colore{blue}{discriminant} G
     &colore{Sepia}{-- used 0.0310744 seconds}
o5 = 140570811483169199470638017932788358544282187717397844656324826769552160278476332
56406502145120855236676811697488882435760217714078399664105019672381338748228576388801
69042329841357623161361759778624522173244483459194112043602458289220741512289591637737
14466361681597648097658753070739833449997864683601657856
i6 : R := ZZ[t,u][w,x,y,z]; pencil = t*sub(F,R) + u*sub(G,R)
             6           6          4           6      6
o7 = (t + u)w  + (t + u)x  + t*w*x*y  + (t + u)y  + u*z
o7 : ZZ[t, u][w, x, y, z]
i8 : time D = &colore{blue}{discriminant} pencil
     &colore{Sepia}{-- used 7.05101 seconds}
                                   375 125                                 374 126 ...
     11918167904272470982401...000t   u    + 44811489377450403137211...000t   u    ...     
o8 : ZZ[t, u]     
i9 : factor D
        125       195    3      2         2      3 30    3      2         2      3 30
o9 = (u)   (t + u)   (25t  + 81t u + 81t*u  + 27u )  (29t  + 81t u + 81t*u  + 27u )  (
18453098603344854356045130076201433820906084922117987408631404035314583354936784858690
19666668055428407222803144055042891867966935429959336227999512218285981355846846846364
626801397625813957058058834010980828766582924640256)
o9 : Expression of class Product
\end{Verbatim}
} \noindent 
In particular, we deduce that the pencil $\langle F,G \rangle$ intersects the discriminant hypersurface 
in $F$ with multiplicity $125$, in $F-G$ with multiplicity $195$, and in other six
 distinct points with multiplicity $30$.
\end{example}
The package also provides methods for working with 
Chow forms and more generally tangential Chow forms 
 of projective varieties (see \cite[p.~104]{gelfandbook} and \cite{GreenMorrison}).
In the following example, we apply some of these methods. 
\begin{example}
Take $C\subset\PP^3$ to be the twisted cubic curve.
{\footnotesize 
\begin{Verbatim}[commandchars=&!$]
i10 : C = kernel &colore!blue$!veronese$(1,3)
              2                       2
o10 = ideal (x  - x x , x x  - x x , x  - x x )
              2    1 3   1 2    0 3   1    0 2
o10 : Ideal of QQ[x , x , x , x ]
                   0   1   2   3
\end{Verbatim}
} \noindent 
The Chow form of $C$ in $\GG(1,3)$ can be obtained as follows:
{\footnotesize 
\begin{Verbatim}[commandchars=&\{\}]
i11 : w = &colore{blue}{chowForm} C
       3                         2      2
o11 = x    - x   x   x    + x   x    + x   x    - 2x   x   x    - x   x   x
       1,2    0,2 1,2 1,3    0,1 1,3    0,2 2,3     0,1 1,2 2,3    0,1 0,3 2,3
      QQ[x   , x   , x   , x   , x   , x   ]
          0,1   0,2   1,2   0,3   1,3   2,3
o11 : --------------------------------------
          x   x    - x   x    + x   x
           1,2 0,3    0,2 1,3    0,1 2,3
\end{Verbatim}
} \noindent 
We can recover $C$ from its Chow form by taking the so-called Chow equations (\cite[p.~102]{gelfandbook}, \cite{cataneseJ}).
{\footnotesize 
\begin{Verbatim}[commandchars=&\{\}]
i12 : C == saturate &colore{blue}{chowEquations} w
o12 = true
\end{Verbatim}
} \noindent 
The $X$-resultant of $C$ can be obtained applying first 
 the duality isomorphism $\GG(1,\PP^3) = \GG(1,{\PP^3}^{\ast})$
and then passing from the Pl\"ucker to the Stiefel coordinates. 
{\footnotesize 
\begin{Verbatim}[commandchars=&\{\}]
i13 : w' = &colore{blue}{dualize} w
       3                         2      2
o13 = x    - x   x   x    + x   x    + x   x    - x   x   x    - 2x   x   x
       0,3    0,2 0,3 1,3    0,1 1,3    0,2 2,3    0,1 1,2 2,3     0,1 0,3 2,3
      QQ[x   , x   , x   , x   , x   , x   ]
          0,1   0,2   1,2   0,3   1,3   2,3
o13 : --------------------------------------
          x   x    - x   x    + x   x
           1,2 0,3    0,2 1,3    0,1 2,3
i14 : &colore{blue}{fromPluckerToStiefel} w' 
         3   3          2   2              2       2          2   3     ...
o14 = - x   x    + x   x   x   x    - x   x   x   x    + x   x   x    - ...
         0,3 1,0    0,2 0,3 1,0 1,1    0,1 0,3 1,0 1,1    0,0 0,3 1,1   ... 
o14 : QQ[x   , x   , x   , x   , x   , x   , x   , x   ]
          0,0   0,1   0,2   0,3   1,0   1,1   1,2   1,3
\end{Verbatim}
} \noindent 
The method \texttt{cayleyTrick} returns a pair 
consisting of the defining ideal of $\PP^1 \times C\subset\PP^7\simeq\PP(\mathrm{Mat}(2,4))$ 
and the  $X$-resultant of $C$, considered as a hypersurface $Z\subset \PP(\mathrm{Mat}(2,4))$.
Theorem~\ref{CayleyTrickTeo} ensures that $Z$ is the dual variety of $\PP^1\times C$.
We can check this using the method \texttt{dualVariety}.
{\footnotesize 
\begin{Verbatim}[commandchars=&\{\}]
i15 : (P1xC,Z) = &colore{blue}{cayleyTrick} C;
i16 : &colore{blue}{dualVariety}(P1xC) == Z
o16 = true
\end{Verbatim}
} \noindent 
\end{example}

\subsection*{Some overlapping packages}
There are two further packages related to resultant computations, which are included in Macaulay2:
\href{https://github.com/Macaulay2/M2/blob/master/M2/Macaulay2/packages/Elimination.m2}{\emph{Elimination}}, 
by M. E. Stillman (see \cite{Elimination}); and 
\href{https://github.com/Macaulay2/M2/blob/master/M2/Macaulay2/packages/EliminationMatrices.m2}{\emph{EliminationMatrices}},
by
N. Botbol, 
L. Bus\'e and
M. Dubinsky (see \cite{EliminationMatrices}).
The former 
contains functions to compute Sylvester resultants. The latter 
can compute different resultant matrices; in particular, 
it contains an implementation of the Macaulay formula.

A further package for working with Chow forms is
 \href{http://page.math.tu-berlin.de/~kohn/packages/Coisotropy.m2}{\emph{Coisotropy}}, 
 by K. Kohn (see \cite{KohnCoisotropic}), which in particular, 
 contains a useful function to compute the degrees 
 of all tangential Chow forms of a given projective variety.
 
 For all these overlapping functions,
 it does not seem easy to rank 
 implementations in terms of efficiency because 
 this generally depends on the problem.
 They also differ 
 in how they handle input and output.
 For instance, the discriminant of a binary form computed using the package \emph{Elimination} 
 lies again in the same ring, rather than in the ring of coefficients,
 and the Chow form of a projective variety computed using \emph{Coisotropy} 
 lies in a polynomial ring, rather than in a quotient ring.
 

\providecommand{\bysame}{\leavevmode\hbox to3em{\hrulefill}\thinspace}
\providecommand{\MR}{\relax\ifhmode\unskip\space\fi MR }
\providecommand{\MRhref}[2]{%
  \href{http://www.ams.org/mathscinet-getitem?mr=#1}{#2}
}
\providecommand{\href}[2]{#2}

\end{document}